\documentclass[reqno,a4paper]{amsart}

\usepackage{mathrsfs,amsmath,amssymb,graphicx}
\usepackage{stmaryrd}
\usepackage[unicode]{hyperref}
\usepackage{latexsym}
\usepackage{layout}
\usepackage[english]{babel}
\usepackage{color}
\usepackage{fancyvrb}
\usepackage{color}
\usepackage{amsmath,amscd}
\usepackage{bm}

\usepackage{subcaption}
\usepackage{xifthen}
\usepackage{transparent}

\usepackage{tikz}
\usetikzlibrary{calc}
\usetikzlibrary{matrix,arrows,decorations}

\theoremstyle{plain}
\newtheorem*{theorem*}{Theorem} 
\newtheorem{theorem}{Theorem}[section]
\newtheorem{lemma}[theorem]{Lemma}

\newtheorem{corollary}[theorem]{Corollary}

\theoremstyle{definition}
\newtheorem{definition}{Definition}[section]
\theoremstyle{remark}
\newtheorem{remark}{Remark}[section]

\newcommand{\HK}{{\rm HK}}

\newcommand{\Sbb} {\mathbb S}

\newcommand{\sphere} {\Sbb^3}
\newcommand{\pair}{(\sphere,\HK)}

\newcommand{\id}{\rm id}
\newcommand{\rel}{{\rm rel\,}}
\newcommand{\PGL}{{\rm PGL}}
\newcommand{\PSL}{{\rm PSL}}

\newcommand{\Compl}[1]{E(#1)}

\newcommand{\Stwo}{\mathfrak{S}^2}

\newcommand{\cout}[1]   {}

\newcommand{\Sym}[2][\sphere]{\mathcal MCG(#1, #2)}
\newcommand{\nbhd}[2]{\mathfrak N(#1; #2)}
\newcommand{\rnbhd}[1]{\mathfrak N(#1)} 
\newcommand{\opennbhd}[2]{\mathring{\mathfrak N}(#1; #2)}
\newcommand{\openrnbhd}[1]{\mathring{\mathfrak N}(#1)}

\newcommand{\pSym}[2][\sphere]{\mathcal MCG_+(#1, #2)}
\newcommand{\MCG}[1]{\mathcal MCG(#1)}
\newcommand{\pMCG}[1]{\mathcal MCG_+(#1)}
\newcommand{\Aut}[1]{\mathcal Homeo(#1)}
\newcommand{\pAut}[1]{\mathcal Homeo_+(#1)}
\newcommand{\Emb}[2]{\mathcal Emb(#1,#2)}

\newcommand{\TSG}[1]{{\rm TSG}(\sphere, #1)}

\definecolor{mygray}{rgb}{0.92,0.92,0.92}

\numberwithin{equation}{section}
\numberwithin{figure}{section}

\title[Unknotting annuli]{Unknotting annuli and handlebody-knot symmetry}
\author{Yi-Sheng Wang}
\address{National Center of Theoretical Sciences, Mathematics Division, Taipei City 106, Taiwan}
\email{yisheng@ncts.ntu.edu.tw}

\date{\today}

\begin{document}

\thanks{}

\begin{abstract}
By Thurston's hyperbolization theorem, 
irreducible handlebody-knots are classified
into three classes: hyperbolic, toroidal, and atoroidal cylindrical. 
It is known that a non-trivial handlebody-knot of genus two 
has a finite 
symmetry group if and only if
it is atoroidal. 
The paper investigates the topology of 
cylindrical handlebody-knots of genus two that
admit an unknotting annulus; we show that 
the symmetry group is trivial if the unknotting annulus is unique and of type $2$. 
\end{abstract}

\maketitle
 
\section{Introduction}\label{sec:intro}
Given a subspace $X$ of an oriented manifold $M$, 
denoted by $(M,X)$, 
its mapping class group $\MCG{M,X}$ is defined as 
the group of 
isotopy classes of self-homeomorphisms of $M$ 
that preserve $X$ setwise, and its 
positive mapping class group $\pMCG{M,X}$ is
the subgroup consisting of 
orientation-preserving homeomorphisms.
When $M=\sphere$ is a $3$-sphere, $\Sym X$ (resp.\ $\pSym X$) 
is referred to as the (resp.\ positive) symmetry group of $X$.

The case where $(\sphere,K)$ is a knot  
has been studied 
by several authors 
(see 
Boileau-Zimmermann \cite{BoiZim:87}, 
Kodama-Sakuma\cite{KodSak:92}, Kawauchi \cite{Kaw:96}, for instance),
and symmetry groups of a large class of knots are now determined.
In particular, if the exterior of a knot is atoroidal,
the symmetry group is finite, and furthermore it is 
either cyclic or dihedral. 
The present paper concerns symmetry groups of handlebody-knots 
of \emph{genus two}, abbreviated to handlebody-knots hereafter;
a knot can be viewed as a handlebody-knot of genus one.  
%
%
We call a handlebody-knot $\pair$ atoroidal, acylindrical or
irreducible\footnote{ 
The standard definition is
that $\pair$ is irreducible 
if there is no $2$-sphere $\Stwo\subset \sphere$
such that $\Stwo\cap\HK$ is an essential disk of $\HK$.
In the genus two case,  irreducibility of $\pair$ 
is equivalent to $\partial$-irreducibility of $\Compl\HK$
by Tsukui \cite[Theorem $1$]{Tsu:75}.
}
if the exterior $\Compl\HK:=\overline{\sphere-\HK}$ 
contains no essential torus, annulus or disks, respectively.
A handlebody-knot $\pair$ is \emph{hyperbolic} if
$\Compl\HK$
admits a complete hyperbolic structure with geodesic boundary. 
By Thurston's hyperbolization theorem, equivariant torus theorem
\cite{Hol:91} and Tollefson's fixed point theorem \cite{Tol:81},
a handlebody-knot is hyperbolic if and only if 
it is irreducible, atoroidal, and acylindrical.

A reducible or toroidal handlebody-knot $\pair$
has an infinite symmetry group. In the case where $\pair$ is trivial, 
Akbas\cite{Akb:08}, Cho \cite{Cho:08} prove that 
the symmetry group, equivalent to the genus two Goeritz group of $\sphere$, 
is finitely presented (see also 
Goeritz \cite{Goe:33}, Scharlemann \cite{Sch:04}); 
finite presentation of symmetry groups of 
non-trivial reducible handlebody-knots
has also been obtained by Koda \cite{Kod:15}.

On the other hand, by the Mostow rigidity theorem, 
the symmetry group of a hyperbolic handlebody-knot is always finite. 
The finiteness of symmetry groups of 
atoroidal \emph{cylindrical} handlebody-knot is recently proved
by Funayoshi-Koda \cite{FunKod:20}. 
Thus a handlebody-knot has a finite symmetry group
if and only if it is non-trivial and atoroidal
\footnote{In the genus two case, non-triviality and atoroidality imply irreducibility.}. Contrary to the case of knots however, less
is known about the structure of these finite symmetry groups.

The present work investigates the structure
of symmetry groups of atoroidal \emph{cylindrical} handlebody-knots.
To state the results, we recall the classification
of essential annuli  
in the exterior $\Compl\HK$ 
of an irreducible atoroidal 
handlebody-knot $\pair$ by Koda-Ozawa \cite{KodOzaGor:15}.
Such annuli are classified into 
four types by \cite[Corollary $3.18$]{KodOzaGor:15}, 
and the classification 
can be described in terms of the boundary of annuli   
\cite[the proof of Theorem $3.3$]{KodOzaGor:15}:
Let $A$ be an essential annulus in $\Compl\HK$.
$A$ is of type $2$ if exactly one component of $\partial A$
bounds a disk in $\HK$. $A$ is of type $3$ if 
no component of $\partial A$ bounds a disk in $\HK$
and there exists an essential disk in $\HK$ 
disjoint from $\partial A$; a type $3$ annulus
can be further classified into two subtypes:
$A$ is of type $3$-$2$ if the components of $\partial A$ are parallel,  
and is of type $3$-$3$ otherwise.
$A$ is of type $4$ if $\partial A$ 
does not bound disks in $\HK$ and no essential disks
in $\HK$ disjoint from $\partial A$.

Recall from Koda \cite{Kod:15} 
an \emph{unknotting} annulus $A$ of $(\sphere,\HK)$ is an
annulus $A\subset\Compl\HK$ such that 
$(\sphere,\HK_A)$ is a trivial handlebody-knot, 
where $\HK_A:=\HK\cup\rnbhd{A}$ and $\mathfrak{N}(A)$
is a regular neighborhood of $A$ in $\Compl\HK$; 
in other words,
$\partial \HK_A$ induces a genus two Heegaard splitting of $\sphere$. 
The existence of unknotting annuli imposes 
topological constraints on $\pair$; 
such constraints are investigated in Section \ref{sec:unknotting}, and
the results therein are summarized in the following. 
\begin{theorem}\label{intro:topology_HK_with_unknotting}
If $A$ is an unknotting annulus of an irreducible handlebody-knot $(\sphere,\HK)$, 
then $A$ is essential, and is either of type $2$ or of type $3$-$3$; 
furthermore, $\Compl\HK$ is atoroidal.
\end{theorem}

Our next result concerns the symmetry group 
of a handlbody-knots $\pair$ admitting
a \emph{unique} unknotting annulus $A$; that is, 
any other unknotting annulus of the same type is 
isotopic to $A$ in $\Compl\HK$. Examples of handlebody knots admitting
a unique unknotting annulus include infinite families of handlebody-knots 
in Motto \cite{Mott:90}, Lee-Lee \cite{LeeLee:12}, and 
Koda \cite[Example $3$]{Kod:15}. 

\begin{theorem}\label{intro:no_symmetries}
An irreducible handlebody-knot admitting 
a unique unknotting annulus of type $2$ has a trivial symmetry group. 
\end{theorem}  

 
\begin{corollary}
An irreducible handlebody-knot admitting a unique unknotting annulus of type $2$ is chiral.
\end{corollary}
Detecting chirality of handlebody knots
is in general a challenging task; several methods are employed
(e.g.\ Motte \cite{Mott:90}, Lee-Lee \cite{LeeLee:12}, Ishii-Iwakiri \cite{IshIwa:12}, Ishii-Iwakiri-Jang-Oshiro \cite{IshIwaJanOsh:13},
Ishii-Kishimoto-Ozawa \cite{IshKisOza:15}) to determine the chirality of 
handlebody knots in 
the Ishii-Kishimoto-Moriuchi-Suzuki knot table \cite{IshKisMorSuz:12}.

Theorem \ref{intro:topology_HK_with_unknotting} can be obtained by standard
$3$-manifold techniques, and classification of knot tunnels of a trefoil knot, whereas Theorem \ref{intro:no_symmetries} 
relies on results on knot-tunnels-preserving 
homeomorphisms by Cho-McCullough \cite{ChoMcC:09}, and 
finiteness theorems on symmetry group
of spatial graphs and handlebody-knots in Cho-Koda \cite{ChoKod:13}
and Funayoshi-Koda \cite{FunKod:20}. 
We also make use of mapping class groups of surfaces 
(see Farb-Margalit \cite{FarMar:11}, \"{O}zg\"{u}r-\c{S}ahin \cite{OzgSah:03}), and 
homotopy types of embedding spaces of subpolyhedra in a surface 
as discussed in Yagasaki \cite{Yag:00}, \cite{Yag:05}.


\section{Preliminaries}\label{sec:prelim}
In this section we fix the convention, and recall some results in
\cite{ChoKod:13}, \cite{ChoMcC:09}, \cite{Kod:15}, \cite{Yag:00},
\cite{Yag:05} needed in subsequent sections.

Throughout the paper, we work in the piecewise linear category,
given subpolyhedra $X_1,\cdots,X_n$ of a manifold $M$, the space 
of self-homeomorphisms of $M$ preserving $X_i$, $i=1,\cdots, n$, 
setwise (resp.\ pointwise) is denoted by 
\[\Aut{M,X_1,\cdots,X_n} \qquad\text{(resp.\ $\Aut{M,\rel X_1,\cdots ,X_n}$ )},\]
and the mapping class group of $(M,X_1,\cdots,X_n)$ is defined as 
\begin{multline*}
\MCG{M,X_1,\cdots,X_n}:=\pi_0(\Aut{M,X_1,\cdots ,X_n})\\
\text{(resp.\ $\MCG{M,\rel X_1,\cdots,X_n}:=\pi_0(\Aut{M,\rel X_1,\cdots ,X_n})$ )}.
\end{multline*} 
The ``+'' subscript is added when considering 
the subspace or the subgroup consisting of orientation-preserving
homeomorphisms
\begin{align*}
\pAut{M,X_1,\cdots, X_n}&\quad \text{(resp.\ $\pAut{M,\rel X_1,\cdots, X_n}$ )}\\
\pMCG{M, X_1,\cdots,X_n}&\quad\text{(resp.\ $\pMCG{M,\rel X_1,\cdots , X_n}$)}.
\end{align*}
If $f$ is a self-homeomorphism of $(M,X_1,\cdots,X_n)$, 
we denote by $[f]$ the mapping class it represents.


Given a subpolyhedron $X$ of $M$, 
$\mathring{X}$ denotes the interior of $X$, and
$\nbhd{X}{M}$ a regular neighborhood of $X$ in $M$, or simply
$\mathfrak{N}(X)$ when $M$ is clear from the context. 
The exterior $\Compl X$ of $X$ in $M$ is the the complement $\opennbhd{X}{M}$
if $X$ has codimension greater than zero, and  
is the closure of $M-X$ otherwise. 
Submanifolds of a manifold $M$ are understood to be proper 
except in some obvious case
where submanifolds are in $\partial M$, and intersection of two submanifolds
are assumed to be transverse.  
A surface in a three-manifold is essential if 
it is non-boundary parallel, incompressible, and $\partial$-incompressible.

\begin{lemma}[\textbf{\cite[Corollary $11.3$, Theorem $16.2$]{ChoMcC:09}}]\label{lm:asymmetry_unknotted_spatial_graph}
Let $\Gamma$ be the union of a tunnel number one knot $K$
or link $L$ and a tunnel $\tau$,
and $f\in\Aut{\sphere,\Gamma}$. Then
\begin{itemize}
\item If $f\in\pAut{\sphere,\Gamma}$ swaps two arcs of $K$, then
$K$ is trivial or $\tau$ is the upper or lower tunnel of a two-bridge knot $K$.

\item If $f\notin\pAut{\sphere,\Gamma}$, then $\tau$ is the tunnel
of either a trivial knot $K$, or a trivial link or Hopf link $L$.
\end{itemize} 
\end{lemma}

Assume $X$ is a subpolyhedron of a surface $\Sigma$ of genus two.
We denote by $\Emb{X}{\Sigma}_0$ 
(resp.\ $\Aut{X}_0$) 
the component of the space of embeddings
of $X$ in $M$ 
(resp.\ self-homeomorphisms of $X$) 
containing the inclusion (resp.\ identity). Let 
\[r:\pAut{\Sigma,X}\rightarrow \pAut{X} 
\quad\text{and}\quad r:\pAut{\Sigma}\rightarrow \Emb{X}{\Sigma}\]
be the restriction maps, which are Serre fibrations
\cite{Hud:66}, \cite{Ham:76}.

\begin{lemma}[\textbf{\cite{Yag:05},\cite{Yag:00}}]\label{lm:embeddings_in_surfaces}
If $\pi_1(X)$ is non-cyclic, then $\Emb{X}{\Sigma}_0$ is contractible.
If $X$ is a finite union of circles, then the natural homomorphism
\[\pi_0\big(r^{-1}(\pAut{X}_0)\big)\rightarrow \pi_0\big(r^{-1}(\Emb{X}{\Sigma}_0)\big)\]
is an isomorphism. 
\end{lemma}
\begin{proof}
The first statement is a special case in \cite[Theorem $1.2$]{Yag:05}.
To see the second statement, we note that by
\cite[the proof of Theorem $1.2$]{Yag:05} and 
the Serre fibration \cite{Yag:00}
\[\Emb{S^1}{\Sigma}_0\rightarrow 
\Emb{X}{\Sigma}_0\rightarrow 
\Emb{X'}{\Sigma}_0\]
where $X=S^1\cup X'$,
the natural map 
\[\pAut{X}_0\rightarrow \Emb{X}{\Sigma}_0\]
is a homotopy equivalence. 
The assertion then follows from the map of Serre fibrations.
\begin{center} 
\begin{tikzpicture}
\node (F1) at (0,2) {$\pAut{\Sigma,\rel X}$};
\node (E1) at (4,2) {$r^{-1}(\pAut{\Sigma}_0)$};
\node (B1) at (8,2) {$\pAut{X}_0$};
\node (F2) at (0,0) {$\pAut{\Sigma,\rel X}$};
\node (E2) at (4,0) {$r^{-1}(\Emb{X}{\Sigma}_0)$};
\node (B2) at (8,0) {$\Emb{X}{\Sigma}_0$};
\draw[->] (F1) -- (E1);
\draw[->] (E1) to (B1);
\draw[->] (F2) to (E2);
\draw[->] (E2) to (B2);
\draw[->] (F1) to node [right]{
\rotatebox[origin=c]{90}{$=$}
} (F2);
 
\draw[->] (E1) to (E2);
\draw[->] (B1) to (B2);
\end{tikzpicture}
\end{center}

\end{proof}

A meridian system $\mathcal{D}$ of a handlebody $H$ 
is a set of disjoint, non-parallel, meridian disks $\{D_1,\cdots,D_n\}$
in $H$ such that the exterior of $\cup_{i=1}^n D_i$
consists of only $3$-balls. Every meridian system
determines a trivalent spine of $H$ \cite{Joh:95}\footnote{The definition here is more restrictive than the one in \cite{Joh:95}.}. 
In particular, given a handlebody-knot $\pair$ and 
a meridian system $\{D_1,D_2,D_3\}$ of $\HK$, 
the induced spine is either a spatial 
$\theta$-curve or handcuff graph.
Given a spatial graph $(\sphere,\Gamma)$, 
$\TSG\Gamma$ denotes the topological symmetry 
group \cite{Sim:86}, which is the image of $\Sym \Gamma$
in $\MCG\Gamma$; note that if $\Gamma$
is a handcuff graph, $\MCG\Gamma$ is the dihedral
group $\mathsf{D}_4$.

The next two lemmas 
follow from the Alexander trick and \cite{Ham:66}, \cite[Section $2$]{Hat:76}, 
\cite[Theorem $1$]{Hat:99}
(see also \cite[Section $2$]{ChoKod:13},\cite[Section $2$]{Kod:15}).

\begin{lemma}\label{lm:symmetry_groups_hk_gamma}
Given a handlebody-knot $\pair$, let $\{D_1,D_2,D_3\}$ be a meridian system of $\HK$, and $\Gamma$ 
the induced spatial graph. 
Then 
\begin{itemize}
\item the natural homomorphism
\[\Sym{\HK,D_1\cup D_2\cup D_3}\rightarrow \Sym\HK\]
is injective;
\item the natural homomorphism give by the Alexander trick
\[\Sym{\HK,D_1\cup D_2\cup D_3}\xrightarrow{\sim} \Sym\Gamma\]
is an isomorphism.
\end{itemize}
\end{lemma}
 
\begin{lemma}\label{lm:A_preserving_homo}
Given a handlebody-knot $\pair$ and an essential annulus $A$ in $\Compl\HK$, 
the natural homomorphism
\[\Sym{\HK,A}\rightarrow \Sym\HK\]
is injective.
\end{lemma} 
 

The next lemma is a direct consequence of 
\cite[Theorems $2.5$ and $3.2$]{Kod:15}.
 
\begin{lemma}\label{lm:spine_atoro_irre_hk}
If $\Gamma$ is a spine of an irreducible atoroidal handlebody-knot 
$\HK$, then $\Sym\Gamma\simeq \TSG\Gamma<\mathsf{D}_4$.
\end{lemma} 
The next theorem, strengthening
the finiteness result in \cite[Theorem $4.4$]{Kod:15}
\footnote{which is sufficient for the present work.},
follows from \cite{FunKod:20}.   
\begin{theorem}[\textbf{\cite{FunKod:20}}]\label{teo:koda_finiteness}
The symmetry group of an atoroidal cylindrical handlebody-knot
is finite.
\end{theorem} 
\begin{remark}
\footnote{The author thanks Yuya Koda for explaining 
to him the work in \cite{FunKod:20} and how
Theorem \ref{teo:koda_finiteness}
is derived therefrom.
}
\cite[the proof of Theorem $4.3$]{FunKod:20} shows that every element in 
the symmetry group of an atoroidal cylindrical handlebody-knot
is of finite order. Theorem \ref{teo:koda_finiteness} then follows from
the injections \cite{ChoKod:13}, \cite{Hat:76}
\[\Sym{\HK}\rightarrow  \MCG{\HK}\rightarrow \MCG{\partial\HK}\]
and $\MCG{\partial\HK}$ being virtually torsion free \cite{BasLub:83}. 
\end{remark} 

We now review some properties of 
the mapping class group of a four-times-punctured sphere $S_{0,4}$.
Up to change of basis, 
the mapping class group $\MCG{S_{0,4}}$
(resp.\ $\pMCG{S_{0,4}}$)
is canonically isomorphic to 
\begin{equation}\label{eq:four_punctures_sphere}
\PGL(2,\mathbb{Z})\ltimes (\mathbb{Z}_2\times\mathbb{Z}_2)
\quad\text{(resp.\ $\PSL(2,\mathbb{Z})\ltimes (\mathbb{Z}_2\times\mathbb{Z}_2)$)}
\end{equation}
via the homomorphisms  
\[
\PGL(2,\mathbb{Z})\rightarrow \MCG{S_{0.4}}
\quad\text{(resp.\ 
$\PSL(2,\mathbb{Z})\rightarrow \pMCG{S_{0,4}}$
)}
\]
given by linear homomorphisms \cite[Section $2.2.5$]{FarMar:11},
where the factor $\mathbb{Z}_2\times \mathbb{Z}_2$ in \eqref{eq:four_punctures_sphere} 
corresponds to the four hyperelliptic involutions.
 
On the other hand, it is known \cite{DrePanShrZha:19}, 
\cite{OzgSah:03} that 
every element of order $2$ in $\PGL(2,\mathbb{Z})$ 
is conjugate to one of the matrices:
\begin{equation}\label{eq:torsion_representatives}
\begin{bmatrix}
-1&0\\
0&1
\end{bmatrix}
\quad
\begin{bmatrix}
0&1\\
1&0
\end{bmatrix}
\quad
\begin{bmatrix}
0&-1\\
1&0
\end{bmatrix}.
\end{equation}
The first two matrices correspond to reflection across
the plane containing all punctures and two punctures, respectively,
while the third matrix corresponds to a rotation by $\pi$ with
two punctures on the axis of rotation.
Similarly, we have every element of order $2$ in $\PSL(2,\mathbb{Z})$
is conjugate to 
$
\begin{bmatrix}
0&-1\\
1&0
\end{bmatrix} 
$ (see \cite[Chapter $7$]{FarMar:11}).
%

Denote by $\sigma$ the composition
\[\MCG{S_{0,4}}\rightarrow \mathsf{S}_4\rightarrow \mathbb{Z}_2\]
induced by even and odd permutations on punctures of $S_{0,4}$,
where $\mathsf{S}_4$ is the permutation group on four punctures.
An element $x\in\MCG{S_{0,4}}$ is \emph{even} if $\sigma(x)=0$,
and is \emph{odd} if $\sigma(x)=1$. 
For instance, 
only the first matrix in \eqref{eq:torsion_representatives} is even.
We call a partition $\mathcal{P}$ of punctures
of $S_{0,4}$ a \emph{grouping} if each part of $\mathcal{P}$
contains two punctures. There are three possible groupings. 
A grouping $\mathcal{P}$
is realized by a loop $l\subset S_{0,4}$ if
$l$ separates one part of $\mathcal{P}$ from the other.

\begin{lemma}\label{lm:even_elements_order_two}
Let $[f]\in\MCG{S_{0,4}}$ be even and of order $2$.
\begin{enumerate}
\item 
If $[f]\in\pMCG{S_{0,4}}$, then $f$ can be isotoped 
such that, for any grouping $\mathcal{P}$, 
$f$ preserves a loop realizing it.
\item
If $[f]\notin\pMCG{S_{0,4}}$, then $f$ can be isotoped 
such that only one grouping cannot be realized by loops 
preserved by $f$.
\end{enumerate}

\end{lemma}
\begin{proof}
If $[f]\in\pMCG{S_{0,4}}$, then $[f]$ is conjugate to
\[
(\id,x)\in
\PSL(2,\mathbb{Z})\ltimes \big(\mathbb{Z}_2\times\mathbb{Z}_2\big),
\]  
namely, a hyperellptic involution, which preserves loops
realizing all three groupings.
If $[f]\notin\pMCG{S_{0,4}}$, then $[f]$ is conjugate to
\[
\big(
\begin{bmatrix}
-1&0\\
0&1
\end{bmatrix},
x\big)
\in
\PGL(2,\mathbb{Z})\ltimes \big(\mathbb{Z}_2\times\mathbb{Z}_2\big),
\]  
that is, the composition of a hyperellptic involution
and a reflection across the plane containing all punctures,
which preserves loops realizing all but one grouping.
\end{proof}


\section{Unknotting annuli}\label{sec:unknotting}



From now on we let $\pair$ be an 
\emph{irreducible} handlebody-knot, and 
$A$ an annulus in $\Compl\HK$.
Recall that $\HK_A=\HK\cup\nbhd{A}{\Compl\HK}$.

\begin{lemma}[\textbf{Essentiality}]\label{lm:essential}
If $\HK_A$ is a handlebody,
then $A$ is essential in $\Compl{\HK}$.
\end{lemma}
\begin{proof}
Observe first that $A$ cannot be boundary-parallel, 
for otherwise, the boundary $\partial\HK_A$ would be non-connected. 
For the same reason, at least one component of $\partial A$
is non-separating and hence essential in $\partial \HK$.
Particularly, if $A$ is compressible, then $\Compl\HK$ 
is $\partial$-reducible, contradicting 
the irreducibility of $\pair$. 
 
Suppose $A$ is $\partial$-compressible, and 
$D$ is a compressing disk. Then the boundary of a regular neighborhood
$\nbhd{A\cup D}{\Compl\HK}$
of $A\cup D$ consists of an annulus parallel to $A$ and a disk $D'$ in
$\Compl{\HK}$.
Since $\pair$ is irreducible, $\partial D'$ bounds a disk in $\partial \HK$.
This implies $A$ is either boundary parallel
or compressible,
but neither can happen.  
\end{proof}

\begin{figure}[h]
\def\svgwidth{0.4\columnwidth}
\begingroup%
  \makeatletter%
  \providecommand\color[2][]{%
    \errmessage{(Inkscape) Color is used for the text in Inkscape, but the package 'color.sty' is not loaded}%
    \renewcommand\color[2][]{}%
  }%
  \providecommand\transparent[1]{%
    \errmessage{(Inkscape) Transparency is used (non-zero) for the text in Inkscape, but the package 'transparent.sty' is not loaded}%
    \renewcommand\transparent[1]{}%
  }%
  \providecommand\rotatebox[2]{#2}%
  \newcommand*\fsize{\dimexpr\f@size pt\relax}%
  \newcommand*\lineheight[1]{\fontsize{\fsize}{#1\fsize}\selectfont}%
  \ifx\svgwidth\undefined%
    \setlength{\unitlength}{1417.32283465bp}%
    \ifx\svgscale\undefined%
      \relax%
    \else%
      \setlength{\unitlength}{\unitlength * \real{\svgscale}}%
    \fi%
  \else%
    \setlength{\unitlength}{\svgwidth}%
  \fi%
  \global\let\svgwidth\undefined%
  \global\let\svgscale\undefined%
  \makeatother%
  \begin{picture}(1,0.596)%
    \lineheight{1}%
    \setlength\tabcolsep{0pt}%
    \put(0,0){\includegraphics[width=\unitlength,page=1]{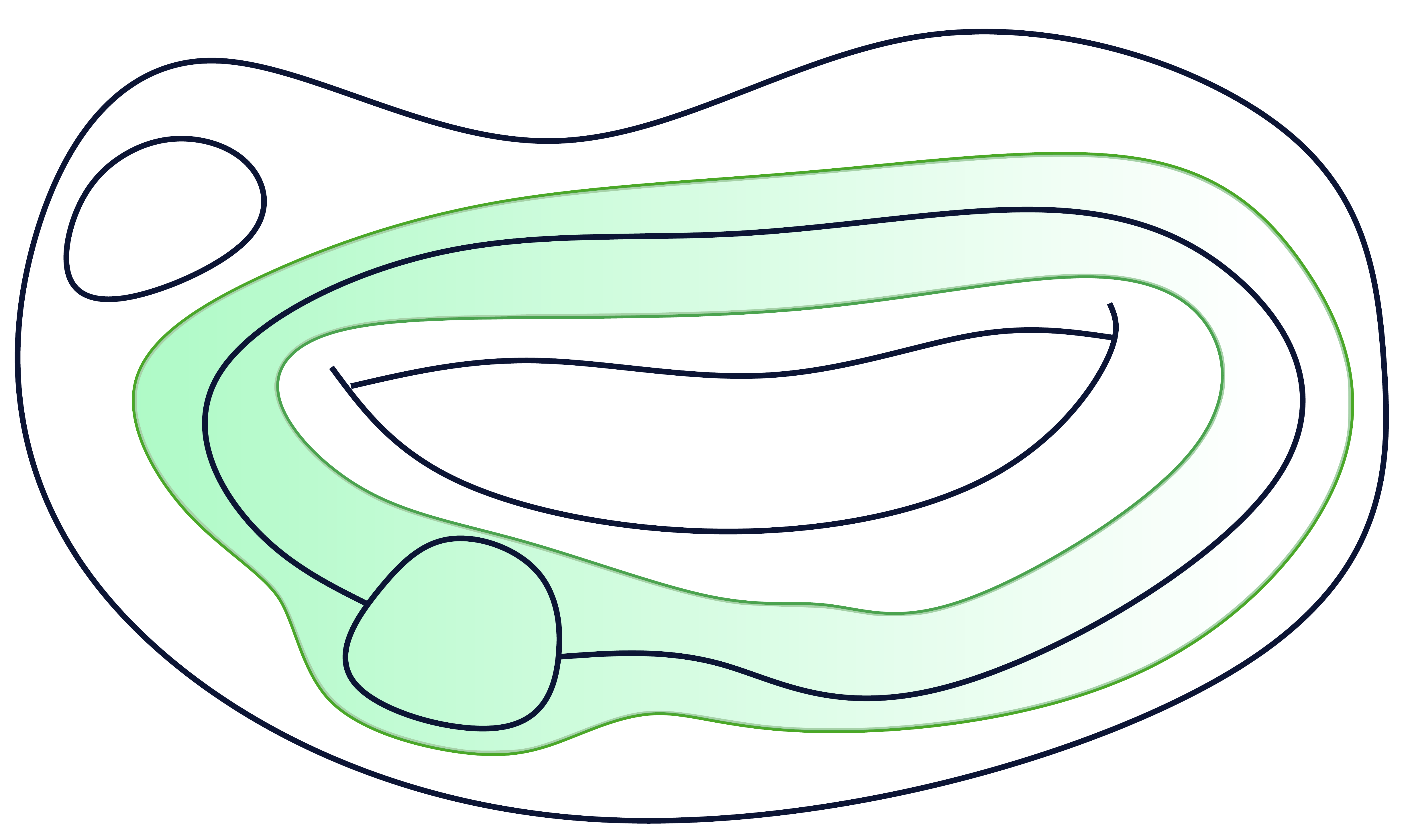}}%
    \put(0.07835863,0.44432921){\color[rgb]{0,0,0}\makebox(0,0)[lt]{\lineheight{1.25}\smash{\begin{tabular}[t]{l}{\tiny $D_{-}$}\end{tabular}}}}%
    \put(0.30816266,0.1171947){\color[rgb]{0,0,0}\makebox(0,0)[lt]{\lineheight{1.25}\smash{\begin{tabular}[t]{l}{\tiny $D_{+}$}\end{tabular}}}}%
    \put(0,0){\includegraphics[width=\unitlength,page=2]{b_compressible.pdf}}%
    \put(0.87637256,0.29427026){\color[rgb]{0,0,0}\makebox(0,0)[lt]{\lineheight{1.25}\smash{\begin{tabular}[t]{l}{\tiny $\alpha$}\end{tabular}}}}%
    \put(0,0){\includegraphics[width=\unitlength,page=3]{b_compressible.pdf}}%
    \put(0.41508517,0.38061099){\color[rgb]{0,0,0}\makebox(0,0)[lt]{\lineheight{1.25}\smash{\begin{tabular}[t]{l}{\tiny $A^\perp$}\end{tabular}}}}%
  \end{picture}%
\endgroup%

\caption{$D_\pm,\alpha$ and $A^\perp$.}
\label{fig:A_perp_b_compressible}
\end{figure}

\begin{lemma}[\textbf{Types of Annuli}]\label{lm:classification}
If $\HK_A$ is a handlebody,
then $A$ is either of type $2$ or of type $3$-$3$.
\end{lemma}
\begin{proof}
Note first that, 
since $\HK_A$ is a handlebody, at most one component of $\partial A$
bounds a disk in $\HK$; also by Lemma \ref{lm:essential}, $A$ is essential.

Let $C\subset A$ be an essential circle in $A$, and denote by 
$A^\perp$ the annulus $C\times I\subset A\times I \simeq \nbhd{A}{\Compl\HK}$,
where $I$ is an interval. $A^\perp$ is a non-separating annulus
in the handlebody $\HK_A$, so it
is compressible, or $\partial$-compressible, or both in $\HK_A$
by \cite[Lemma $9$]{BonOta:83}, \cite[Lemma $2.4$]{HayShi:01}.

If $A^\perp$ is compressible in $\HK_A$, exactly one 
component of $\partial A$ bounds a meridian disk in $\HK$
since $A$ is essential.
This implies $A$ is of type $2$.
If $A^\perp$ is incompressible, then it 
is necessarily $\partial$-compressible. Compressing $A^\perp$
with a compressing disk $D'$, we obtain a non-separating disk $D$
in $\HK_A$ disjoint from $A$.
In particular, 
$A^\perp$ can be obtained by the boundary of a
regular neighborhood of the union of
$D$ and an arc $\alpha\subset\partial \HK_A$ transverse to 
$D'$ and connecting to the same side of $D$.

Since $A^\perp$ is incompressible in $\HK_A$, it is also 
incompressible in the solid torus $V=\HK_A-\opennbhd{D}{\HK_A}$, 
and therefore 
$A^\perp$ separates $V$ into two solid tori $V_+,V_-$ 
with two disk components 
$D_\pm = \partial\mathfrak{N}(D) -\partial \HK$
in $V_\pm$ (Fig.\ \ref{fig:A_perp_b_compressible}) and the core of $A^\perp$
some multiples of longitudes of $V_+$ and $V_-$, respectively.
This implies $A$ is of type $3$-$3$.
\end{proof}
 
\begin{theorem}\label{teo:atoroidality}
If $A$ is unknotting,
then $\pair$ is atoroidal.
\end{theorem}
\begin{proof}
We prove by contradiction. Suppose $T$ 
is an incompressible torus in $\Compl\HK$
such that $\#A\cap T$ is minimized. 
By Lemma \ref{lm:essential}, $A$ is essential,
so every component of $A\cap T$ is essential 
in $A$ and $T$. Since $\HK_A$ is atoroidal
and $T\cap\HK=\emptyset$, $\# T\cap A$
is a positive even number. 
And by Proposition \ref{lm:classification},
$A$ is either of type $2$ or of type $3$-$3$. 

\noindent
\textbf{Case $1$: $A$ is of type $2$.}
Let $B\subset T$ be an annulus with $B\cap A=\partial B$,
and $B'\subset A$ be the annulus cut off by $\partial B$ (Fig.\ \ref{fig:atoroidality_type_two}).
Push $B\cup B'$ slightly away from $\HK$, 
we obtain a torus $T_B\subset\Compl{\HK_A}$.
Let $V$ be the closure of the component of $\sphere-T_B$ 
with $V\cap \HK=\emptyset$. Then $\HK_A$ being atoroidal implies that
$V$ is a solid torus. Since $A$ is of type $2$, 
$\partial B'$ bounds a disk in $\sphere-\mathring{V}$.
Therefore the core of $B'$ is a longitude of $V$,
and hence
\[\pi_1(B')\rightarrow \pi_1(V)\]
is an isomorphism. In particular, $B'$
is parallel to $B$ through $V$, so 
one can isotope $A$ to decrease $\#A\cap T$,
a contradiction.

\begin{figure}[t]
\begin{subfigure}{0.45\textwidth}
\centering
\def\svgwidth{.9\columnwidth}
\begingroup%
  \makeatletter%
  \providecommand\color[2][]{%
    \errmessage{(Inkscape) Color is used for the text in Inkscape, but the package 'color.sty' is not loaded}%
    \renewcommand\color[2][]{}%
  }%
  \providecommand\transparent[1]{%
    \errmessage{(Inkscape) Transparency is used (non-zero) for the text in Inkscape, but the package 'transparent.sty' is not loaded}%
    \renewcommand\transparent[1]{}%
  }%
  \providecommand\rotatebox[2]{#2}%
  \newcommand*\fsize{\dimexpr\f@size pt\relax}%
  \newcommand*\lineheight[1]{\fontsize{\fsize}{#1\fsize}\selectfont}%
  \ifx\svgwidth\undefined%
    \setlength{\unitlength}{1360.62992126bp}%
    \ifx\svgscale\undefined%
      \relax%
    \else%
      \setlength{\unitlength}{\unitlength * \real{\svgscale}}%
    \fi%
  \else%
    \setlength{\unitlength}{\svgwidth}%
  \fi%
  \global\let\svgwidth\undefined%
  \global\let\svgscale\undefined%
  \makeatother%
  \begin{picture}(1,0.58333333)%
    \lineheight{1}%
    \setlength\tabcolsep{0pt}%
    \put(0,0){\includegraphics[width=\unitlength,page=1]{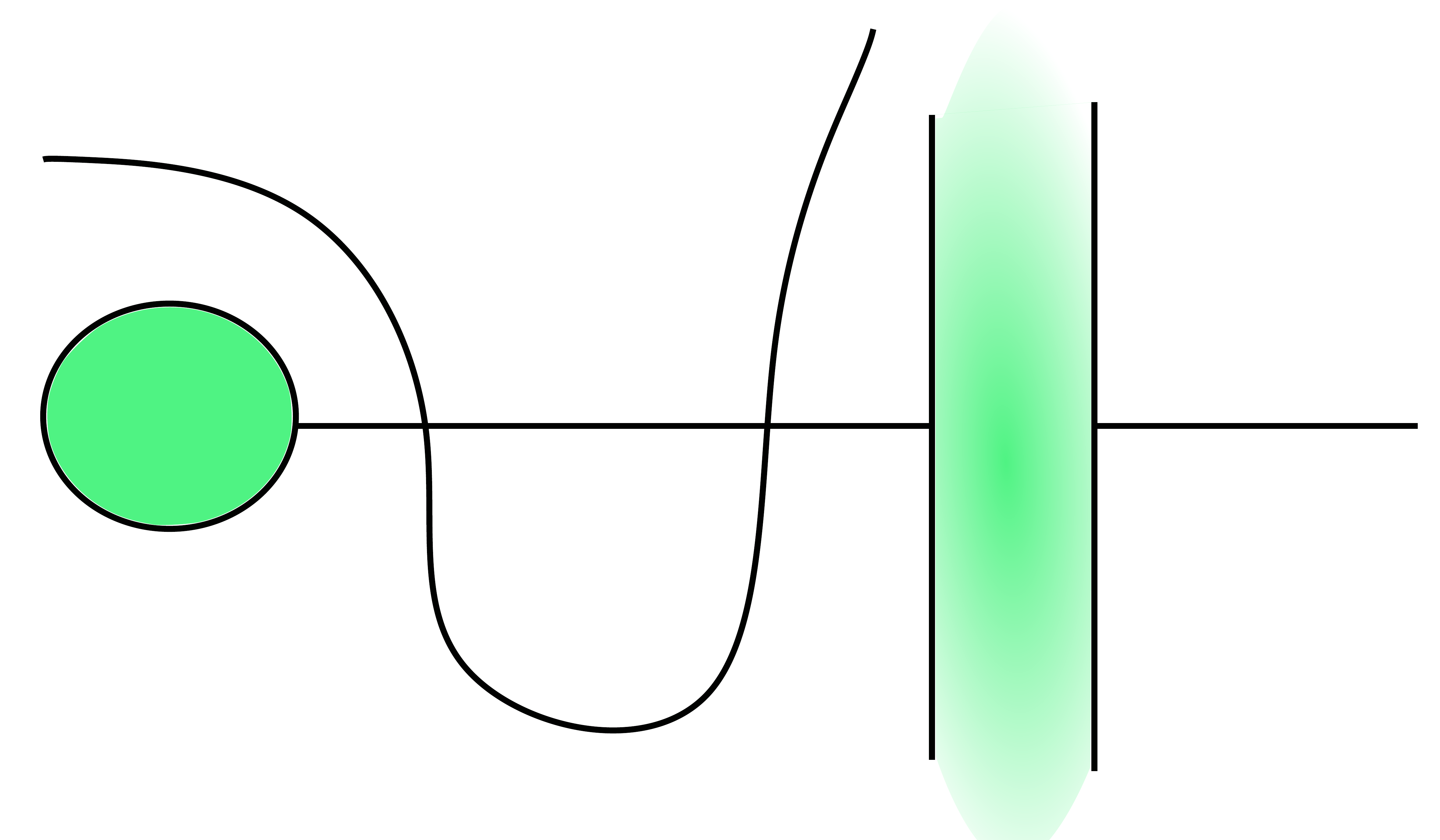}}%
    \put(0.05222892,0.29068718){\color[rgb]{0,0,0}\makebox(0,0)[lt]{\lineheight{1.25}\smash{\begin{tabular}[t]{l}{\footnotesize $\HK$}\end{tabular}}}}%
    \put(0.33026679,0.04503187){\color[rgb]{0,0,0}\makebox(0,0)[lt]{\lineheight{1.25}\smash{\begin{tabular}[t]{l}{\tiny $B$(innermost)}\end{tabular}}}}%
    \put(0.55364634,0.25094452){\color[rgb]{0,0,0}\makebox(0,0)[lt]{\lineheight{1.25}\smash{\begin{tabular}[t]{l}{\tiny $A$}\end{tabular}}}}%
    \put(0.78600216,0.24604575){\color[rgb]{0,0,0}\makebox(0,0)[lt]{\lineheight{1.25}\smash{\begin{tabular}[t]{l}{\tiny $A$}\end{tabular}}}}%
    \put(0.07841008,0.4901228){\color[rgb]{0,0,0}\makebox(0,0)[lt]{\lineheight{1.25}\smash{\begin{tabular}[t]{l}{\tiny $T$}\end{tabular}}}}%
    \put(0.34604034,0.30401907){\color[rgb]{0,0,0}\makebox(0,0)[lt]{\lineheight{1.25}\smash{\begin{tabular}[t]{l}{\tiny $B'$}\end{tabular}}}}%
    \put(0.6584748,0.14055697){\color[rgb]{0,0,0}\makebox(0,0)[lt]{\lineheight{1.25}\smash{\begin{tabular}[t]{l}{\footnotesize $\HK$}\end{tabular}}}}%
  \end{picture}%
\endgroup%

\caption{Type $2$ annulus $A$.}
\label{fig:atoroidality_type_two}
\end{subfigure}
\begin{subfigure}{0.52\textwidth}
\centering
\def\svgwidth{.9 \columnwidth}
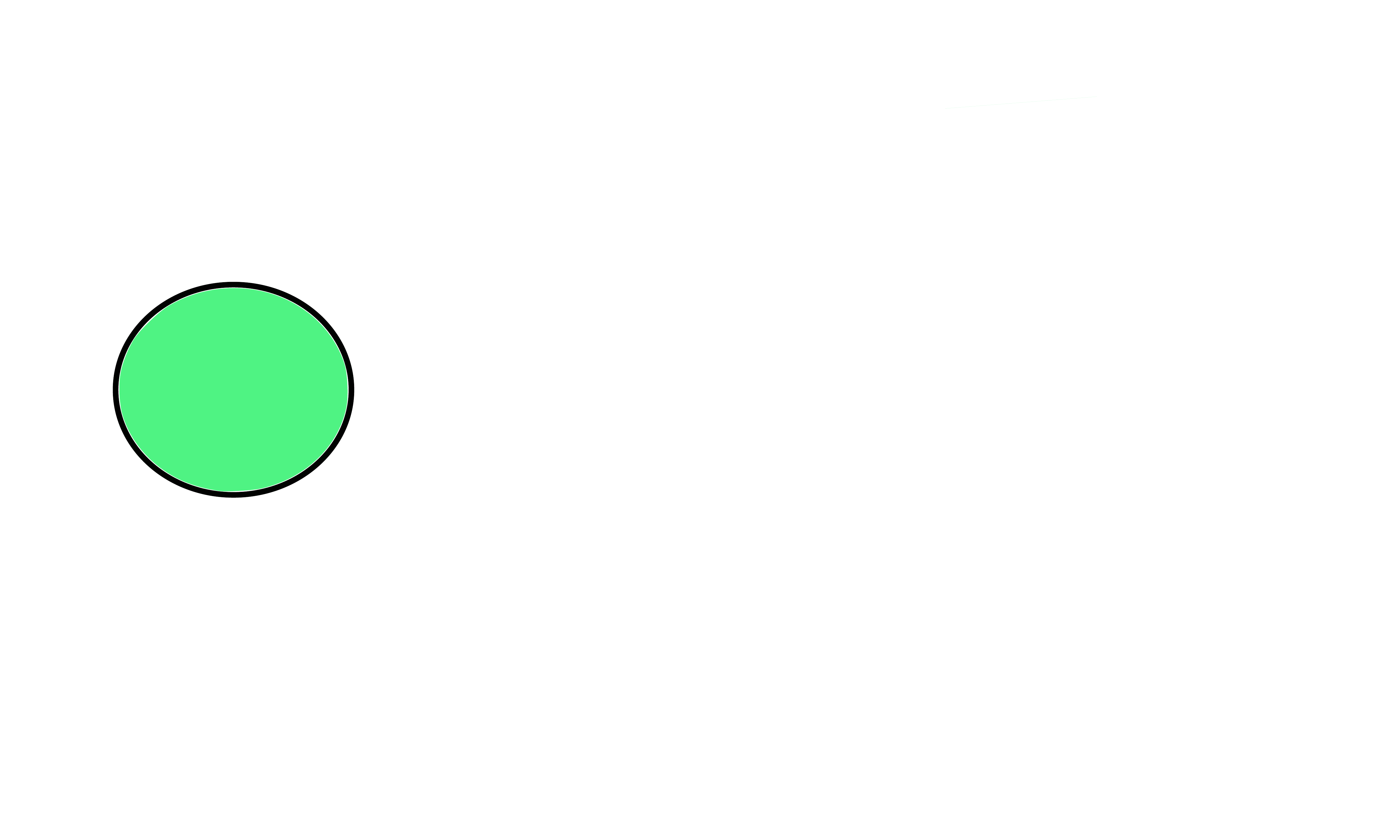
\caption{Type $3$-$3$ annulus $A$ with $A\cap T$ meridional.}
\label{fig:atoroidality_type_threethree_1}
\end{subfigure}
\caption{Schematic diagrams (Types $2$ and $3$-$3$).}
\end{figure}

\noindent
\textbf{Case $2$: $A$ is of type $3$-$3$.}
Let $U\supset \HK$ be the solid torus bounded by $T$.

\textbf{Subcase $2i$: $T\cap A\subset U$ is meridional.} 
There exists an annulus $B\subset A$ with $B\cap T=\partial B$
and $B\not\subset U$. Let $B'\subset T$ be an annulus 
cut off by $\partial B$ (Fig.\ \ref{fig:atoroidality_type_threethree_1}).
Let $T_B$ be the union $B\cup B'$, and 
$V$ be the component of $\sphere-T_B$ not containing $\HK$.
Since $T_B$ has less intersection with $A$ than $T$ does, 
$V$ is necessarily a solid torus.
On the other hand, because 
$T\cap A$ are meridional in $\partial U$, 
the core of $B'$ is a longitude of $V$, and hence $B',B$ are parallel
through $V$. Isotoping $T$ through $V$ gives a contradiction to 
the minimality of $\# A\cap T$.

\begin{figure}[b]
\begin{subfigure}{0.48\textwidth}
\centering
\def\svgwidth{.9\columnwidth}
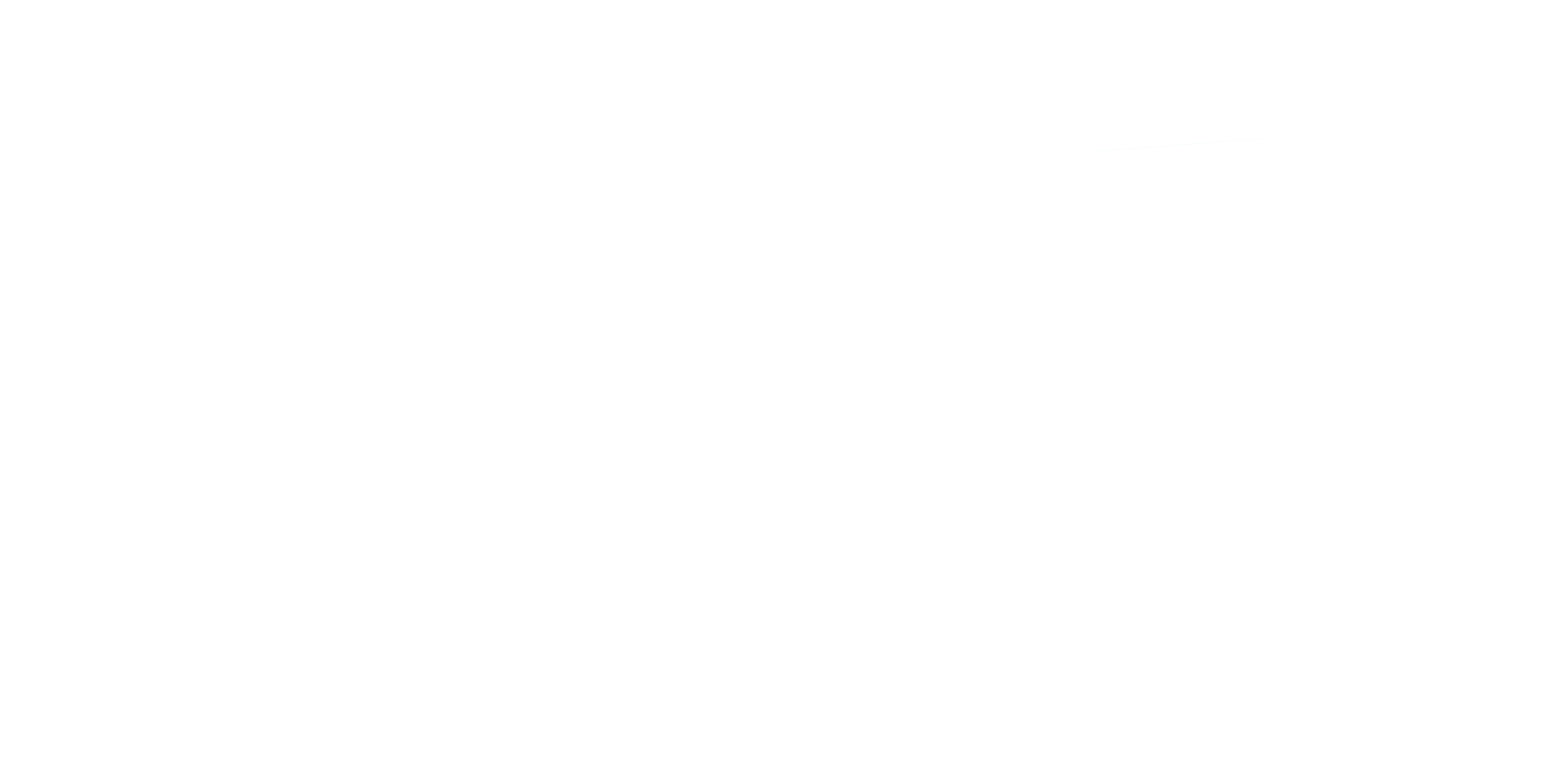
\caption{Type $3$-$3$ annulus $A$ with $\# A\cap T>2$.}
\label{fig:atoroidality_type_threethree_2}
\end{subfigure}
\begin{subfigure}{0.48\textwidth}
\centering
\def\svgwidth{.75 \columnwidth}
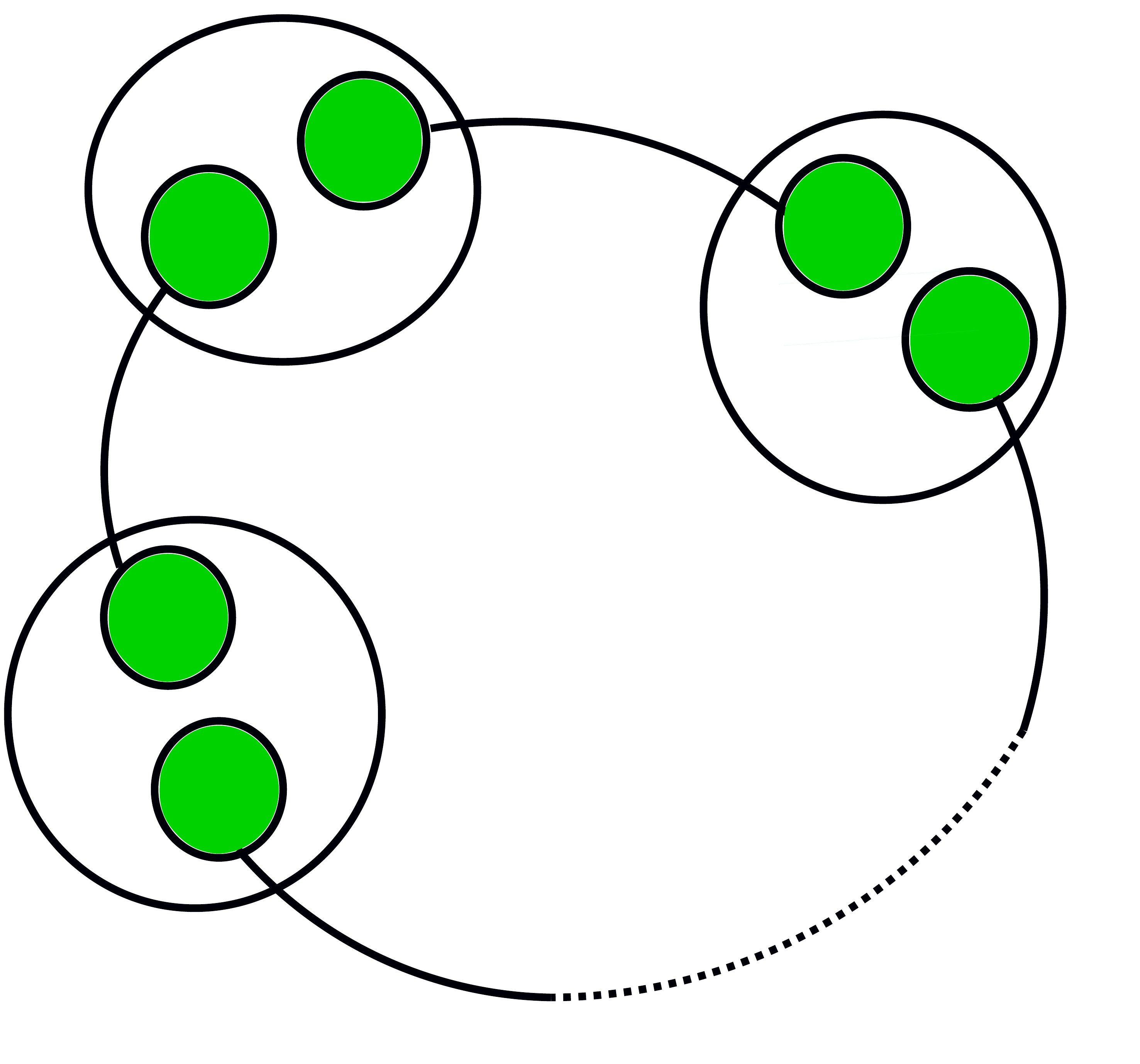
\caption{Type $3$-$3$ annulus $A$ with $\#A\cap T=2$.}
\label{fig:atoroidality_type_threethree_3}
\end{subfigure}
\caption{Schematic diagrams (Type $3$-$3$).}
\end{figure}

\textbf{Subcase $2ii$: $T\cap A\subset U$ is non-meridional.}
We first prove that $\# T\cap A$ is at most $2$. If
$\# T\cap A>2$, then there is an annulus $B\subset A$
with $B\cap T=\partial B$ and $B\subset U$.
Since $\partial B\subset U$ is not meridional,
$B$ divides $U$ into two solid tori $U_1, U_2$ (Fig.\ \ref{fig:atoroidality_type_threethree_2}), 
and the boundary of the one containing $\HK$ is 
an incompressible torus in $\Compl\HK$ 
having 
less intersection with $A$, contradicting the minimality.
%
%

Suppose $\# T\cap A=2$, and hence $T\cap A$ cuts $T$ into two annuli $A_1,A_2$.
Denote by $A_m\subset A$ the annulus with $A_m\cap T=\partial A_m$.
Note that $A_m$ is necessarily in $\Compl U$.
By the atoroidality of $\HK_A$, the components $V_i$ bounded by 
$A_m\cup A_i$, $i=1,2$, with $\HK\not\subset V_i$ are solid tori (Fig.\ \ref{fig:atoroidality_type_threethree_3}).
If one of $A_m\rightarrow V_i$, $i=1$ or $2$, induces an isomorphism
on $\pi_1$. Then $V_1\cup V_2=\Compl U$,
is a solid torus, contradicting the incompressibility of $T=\partial U$.
On the other hand, if neither of $A_m\rightarrow V_i$, $i=1,2$,
induces an isomorphism on $\pi_1$, the core of $U$, and hence the core of $A_m$,
is a torus knot by the classification of Seifert fiber structure  
of $\sphere$ \cite{Sei:33}. Therefore $\HK_A$ can be identified
with $\nbhd{A_m}{\Compl U}\cup\nbhd{\gamma}{U}$, where $\gamma$
is an arc in $U$ with $\partial \gamma\subset\partial A_m$. 
Since $\HK_A$ is trivial,
$\gamma$ can be viewed as a tunnel
of a torus knot (the core of $A_m$).

By the classification of tunnels of a torus knot \cite{BoiRosZie:88}, 
\cite{Mor:88} (see also \cite{ChoMcC:09}), 
there is $f\in\Aut{\sphere,\nbhd{A_m}{\Compl U}}$
sending $f(\gamma)$ to an arc $t$ isotopic to 
an essential arc in $A_1$ (or $A_2$) in $\Aut{U,\rel U}$.
One can further isotope $f$ such that it preserves $A_1$  
with $f$ still sending $\gamma$ to $t$.
Thus $\pair$ can be identified with 
a regular neighborhood of the union of $\partial A_m=\partial A_i$
and an essential arc of $A_i$, $i=1$ or $2$, in $\sphere$, 
which is reducible, however.

\end{proof}




\section{Symmetry groups}\label{sec:symmetries}
 
\begin{definition}
$\pair$ admits 
a unique unknotting annulus $A$ if
given another unknotting annulus $A'$ of the same type, 
there exists a path $F_t\in \Aut{\sphere,\HK}$ 
such that $F_0=\id$ and $F_1(A')=A$.  
\end{definition}  
If $\pair$
admits a unique unknotting annulus $A$, 
there are isomorphisms
\begin{equation}\label{eq:isomorphism_sym_hk_A}
\Sym{\HK,\mathfrak{N}(A)}\simeq\Sym{\HK,A}\xrightarrow{\sim} \Sym\HK
\end{equation}
by Lemma \ref{lm:A_preserving_homo}. Before restating the 
main result, we recall that 
essential annuli of type $2$ can be divided 
into two subtypes \cite{FunKod:20}.
 
\begin{definition}
Let $A$ be a type $2$ unknotting annulus of $\pair$, and
$D\subset \HK$ be a meridian disk bounded by 
a component of $\partial A$. 
Then $A$ is of type $2$-$1$ if $D$ is separating,
and is of type $2$-$2$ otherwise.
\end{definition}

\begin{theorem}\label{teo:no_symmetries}
If $\pair$
admits a unique unknotting annulus of type $2$,
then $\Sym\HK=1$.
\end{theorem}

\begin{proof} 
 
\textbf{Case $1$: $A$ is of type $2$-$1$.}
In this case, $A$ determines a meridian system $\{D_1,D_2,D_3\}$, where
$D_3$ is a separating disk bounding a component $\alpha_3$ 
of $\partial A$, and $D_1,D_2$ are meridian disks 
of the two tori $\HK-\opennbhd{D_3}{\HK}$, respectively, one of which, say $D_1$, 
intersecting with the other component $\alpha_1$ of $\partial A$ at a single point (Fig.\ \ref{fig:typetwoi}).
Since $A$ is unique, $D_1,D_2,D_3$ are unique in $\HK$ up to isotopy, and
therefore the injection 
\begin{equation}\label{eq:isom_dual_disk_hk}
\Sym{\HK,D_1\cup D_2\cup D_3}\rightarrow \Sym\HK
\end{equation} 
is an isomorphism. 
On the other hand, by Lemma \ref{lm:spine_atoro_irre_hk}, if 
$\Gamma$ is a spatial graph associated to the meridian system $(D_1,D_2,D_3)$
of $\HK$, then there are isomorphisms 
\begin{equation}\label{eq:isom_graph_dual_disk}
\TSG\Gamma\simeq 
\Sym\Gamma\simeq \Sym{\HK,D_1\cup D_2\cup D_3}
\end{equation}
\eqref{eq:isom_dual_disk_hk} and \eqref{eq:isom_graph_dual_disk} together 
imply $\Sym\HK<\mathsf{D}_4$. 

\begin{figure}[b]
\begin{subfigure}{0.47\textwidth}
\centering
\def\svgwidth{.8\columnwidth}
\begingroup%
  \makeatletter%
  \providecommand\color[2][]{%
    \errmessage{(Inkscape) Color is used for the text in Inkscape, but the package 'color.sty' is not loaded}%
    \renewcommand\color[2][]{}%
  }%
  \providecommand\transparent[1]{%
    \errmessage{(Inkscape) Transparency is used (non-zero) for the text in Inkscape, but the package 'transparent.sty' is not loaded}%
    \renewcommand\transparent[1]{}%
  }%
  \providecommand\rotatebox[2]{#2}%
  \newcommand*\fsize{\dimexpr\f@size pt\relax}%
  \newcommand*\lineheight[1]{\fontsize{\fsize}{#1\fsize}\selectfont}%
  \ifx\svgwidth\undefined%
    \setlength{\unitlength}{850.39370079bp}%
    \ifx\svgscale\undefined%
      \relax%
    \else%
      \setlength{\unitlength}{\unitlength * \real{\svgscale}}%
    \fi%
  \else%
    \setlength{\unitlength}{\svgwidth}%
  \fi%
  \global\let\svgwidth\undefined%
  \global\let\svgscale\undefined%
  \makeatother%
  \begin{picture}(1,0.83333333)%
    \lineheight{1}%
    \setlength\tabcolsep{0pt}%
    \put(0,0){\includegraphics[width=\unitlength,page=1]{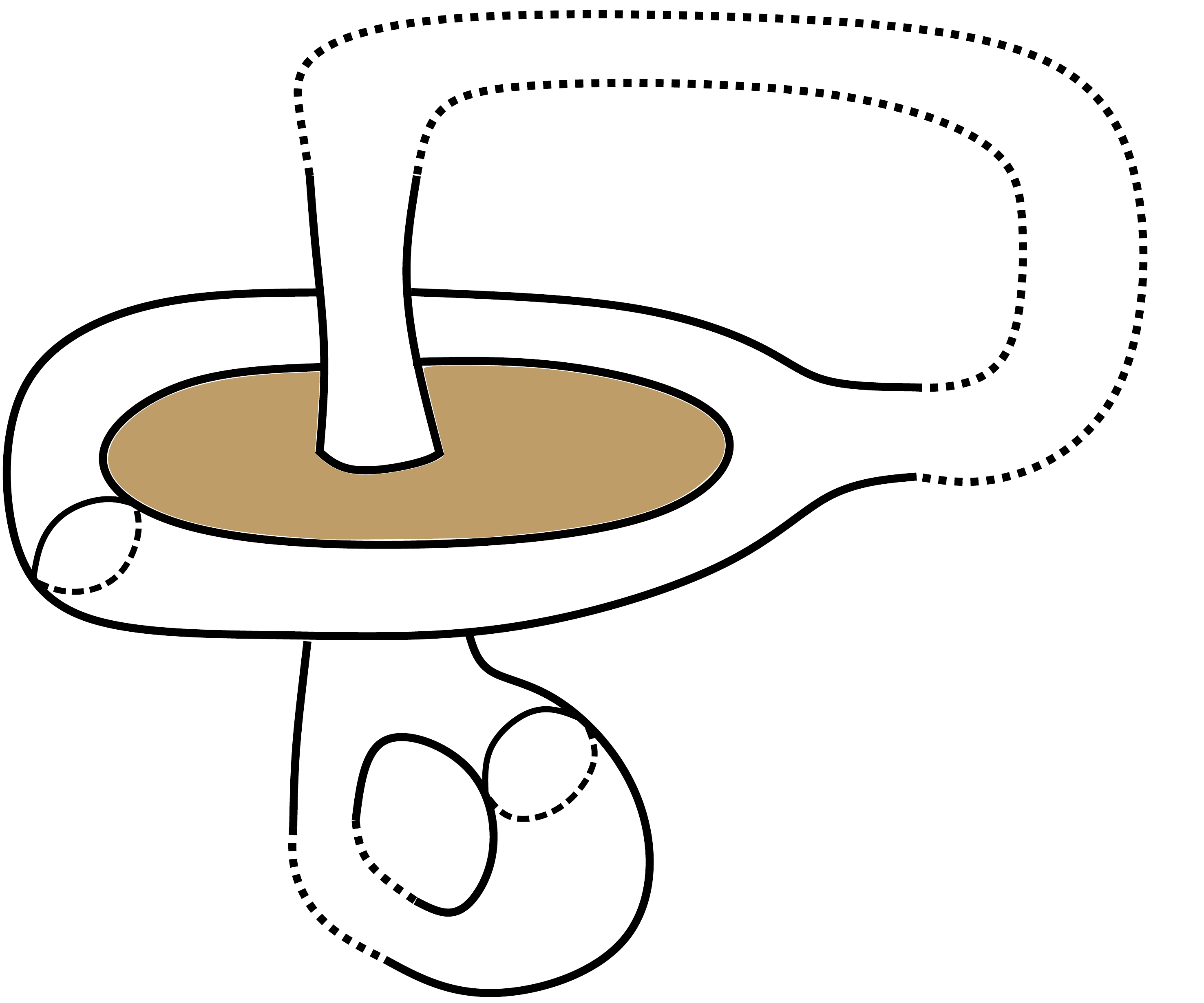}}%
    \put(0.37386507,0.47051919){\color[rgb]{0,0,0}\makebox(0,0)[lt]{\lineheight{1.25}\smash{\begin{tabular}[t]{l}{\tiny $\alpha_3$}\end{tabular}}}}%
    \put(0.3060947,0.3538244){\color[rgb]{0,0,0}\makebox(0,0)[lt]{\lineheight{1.25}\smash{\begin{tabular}[t]{l}{\tiny $\alpha_1$}\end{tabular}}}}%
    \put(0.1832654,0.40782698){\color[rgb]{0,0,0}\makebox(0,0)[lt]{\lineheight{1.25}\smash{\begin{tabular}[t]{l}{\tiny $A$}\end{tabular}}}}%
    \put(0,0){\includegraphics[width=\unitlength,page=2]{typetwoi.pdf}}%
    \put(0.40832466,0.1900551){\color[rgb]{0,0,0}\makebox(0,0)[lt]{\lineheight{1.25}\smash{\begin{tabular}[t]{l}{\tiny $D_2$}\end{tabular}}}}%
    \put(0,0){\includegraphics[width=\unitlength,page=3]{typetwoi.pdf}}%
    \put(0.03799049,0.36941926){\color[rgb]{0,0,0}\makebox(0,0)[lt]{\lineheight{1.25}\smash{\begin{tabular}[t]{l}{\tiny $D_1$}\end{tabular}}}}%
    \put(0,0){\includegraphics[width=\unitlength,page=4]{typetwoi.pdf}}%
  \end{picture}%
\endgroup%

\caption{Type $2$-$1$ annulus.}
\label{fig:typetwoi}
\end{subfigure}
\begin{subfigure}{0.47\textwidth}
\centering
\def\svgwidth{.8\columnwidth}
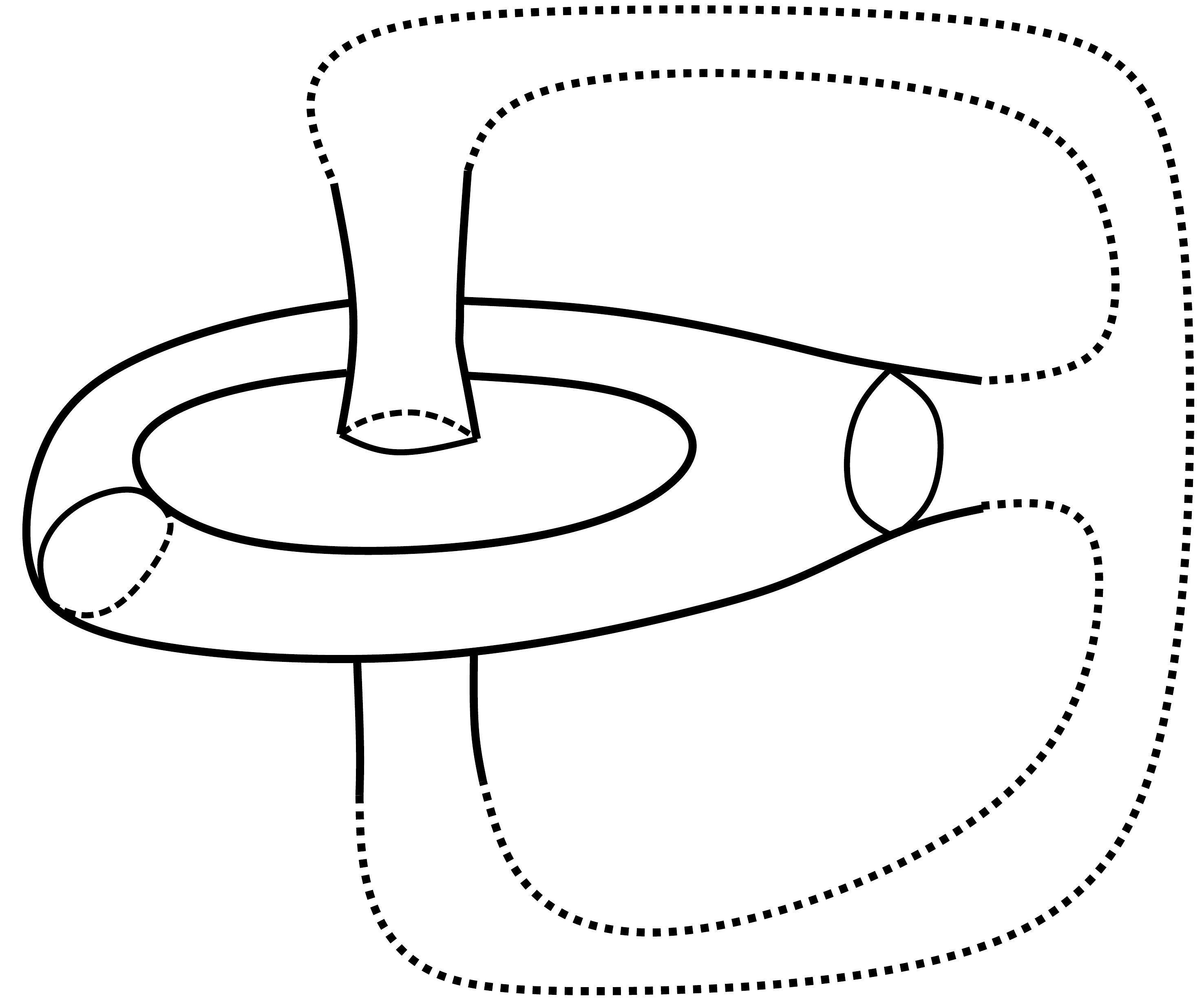
\caption{Type $2$-$2$ annulus.}
\label{fig:typetwoii}
\end{subfigure}
\caption{}
\end{figure}
\centerline{ 
{\bf Claim: $\pSym\HK$ is trivial.} 
}
Given $f\in\pAut{\sphere,\HK}$,
it may be assumed that $f\in\pAut{\sphere,\HK,D_1\cup D_2 \cup D_3}$ by \eqref{eq:isom_dual_disk_hk}, and $f(A)=A$ by the uniqueness of $A$.
$f$ does not permute $\{D_1,D_2,D_3\}$, 
and cannot reverse two sides of $D_1$ and $D_3$ 
since $f$ is orientation-preserving. 
We may hence assume that $f$ 
restricts to the identity on $D_1$ and $\alpha_1$.

On the other hand, $f$ might reverse two sides of $D_2$. By \eqref{eq:isom_dual_disk_hk},
if $f$ fixes two sides of $D_2$, then
$[f]$ is the identity in 
$\pMCG{\sphere,\HK}$,
and if $f$ reverses two sides of $D_2$,
$[f]$ is an element of order $2$ 
in $\Sym\HK$.

Now suppose $[f]$ is non-trivial in 
$\pSym{\HK}$. 
By Theorem \ref{teo:atoroidality} and \cite[Lemma $2.3$]{ChoKod:13},
$\pMCG{\Compl\HK,\partial\Compl\HK}$ is trivial, and therefore 
the composition 
\begin{equation}\label{eq:composition_injections}
\pMCG{\sphere,\HK}\rightarrow \pMCG{\HK}\rightarrow \pMCG{\partial\HK}
\end{equation}
is injective. In particular,
$[f\vert_{\partial \HK}]$ is an element of order $2$
in $\pMCG{\partial\HK}$.
On the other hand,  
there is a Serre fibration  
\begin{multline*}
\pAut{\partial \HK, \rel \alpha_1\cup \partial D_1}\rightarrow 
r^{-1}(\Emb{\alpha_1\cup \partial D_1}{\partial \HK}_0)\\
\xrightarrow{r} \Emb{\alpha_1\cup \partial D_1}{\partial \HK}_0 
\end{multline*}
derived from \cite[Theorem $1.1$]{Yag:00}, where 
$r:\Aut{\partial\HK}\rightarrow \Emb{\alpha_1\cup\partial D_1}{\partial\HK}$
is the restriction map. 
Since $\Emb{\alpha_1\cup \partial D_1}{\partial \HK}_0$
is contractible \cite{Yag:05}, the homomorphism 
\[\pMCG{\partial \HK, \rel \alpha_1\cup\partial D_1}\rightarrow \pMCG{\partial\HK}\]
is injective, and hence 
$[f\vert_{\partial \HK}]$
is an element of order $2$ in 
\[\pMCG{\partial\HK, \rel \alpha_1\cup\partial D_1},\]
which however, is isomorphic to the torsion free group
\[ 
\pMCG{\partial\HK-\openrnbhd{\alpha_1\cup\partial D_1}}.
\]
This gives a contradiction, and therefore $\pSym\HK=1$.
\cout{
by the map of Serre fibrations

\begin{figure}[h]
\begin{tikzpicture}
\node (F1) at (0,0) {$\pAut{\partial\HK,\rel \alpha_1\cup \partial D_1}$};
\node (E1) at (4.5,0) {$\pAut{\partial\HK}_0$}; 
\node (B1) at (8.3,0) {$\Emb{\alpha_1\cup \partial D_1}{\partial\HK}$};
\node (F2) at (0,2) {$\pAut{\partial\HK,\rel \rnbhd{N}(\alpha_1\cup \partial D_1)}$};
\node (E2) at (4.6,2) {$\pAut{\partial\HK}_0$};
\node (B2) at (8.3,2) {$\Emb{\rnbhd{\alpha_1\cup \partial D_1}}{\partial\HK}$};
\draw[->] (F1) -- (E1);
\draw[->] (E1) to (B1);
\draw[->] (F2) to (E2);
\draw[->] (E2) to (B2);
\draw[->] (F2) to (F1);
\draw[->] (E2) to (E1);
\draw[->] (B2) to (B1);
\end{tikzpicture}
\end{figure}
%
%
%
and hence isomorphic to the torsion free group
\[\pMCG{\partial\HK-\openrnbhd{\alpha_1\cup\partial D_1}},\] 
contradicting that $f\vert_{\partial\HK}$ induces an element of order $2$. 
}
\smallskip
 
\centerline{
{\bf Claim: $\Sym\HK=1$.}
}
Suppose $\Sym\HK$ is non-trivial, and $f\in\Aut{\sphere,\HK}$
is an orientation-reversing homeomorphim. 
It may be assumed that $f(\mathfrak{N}(A))=\mathfrak{N}(A)$
by the uniqueness of $A$. Denote by 
$D_3^\pm$ two non-parallel meridian disks in $\HK_A$
that are parallel to $D_3$ in $\HK$. Then $\{D_2,D_3^+,D_3^-\}$
is a meridian disk system of $\HK_A$; let $\Gamma_A$ be
the associated handcuff graph.

Since $f$ preserves $D_3$, one can isotope $f$ such that
it preserves $D_3^\pm$, and hence $f$ 
represents an element in
\[\Sym{\HK_A,D_2, D_3^+, D_3^-}\simeq \Sym{\Gamma_A}.\] 
Notice that $f$ does not permute $\{D_2,D_3^+,D_3^-\}$
since it keeps the two sides of $D_3$.

Let $L$ be constituent link of $\Gamma_A$, and $\gamma$
the connecting arc. Then $\gamma$
is an unknotting tunnel of $L$, and $f$ induces an orientation-reversing 
homeomorphism preserving $L,\gamma$. By
Lemma \ref{lm:asymmetry_unknotted_spatial_graph}, 
$\gamma$ is the tunnel of 
either a trivial link or a Hopf link (Fig.\ \ref{fig:tunnel_hopf}). In the former, 
$\pair$ is trivial, contradicting the assumption, whereas in the latter, 
$\pair$ is the handlebody-knot $4_1$ in \cite[Table $1$]{IshKisMorSuz:12} (Fig.\ \ref{fig:hk4_1}), which has two non-isotopic unknotting annuli of type $2$. 
Therefore $\Sym\HK=1$. 
\begin{figure}[h]
\begin{subfigure}{0.39\textwidth}
\centering
\def\svgwidth{.65\columnwidth}
\begingroup%
  \makeatletter%
  \providecommand\color[2][]{%
    \errmessage{(Inkscape) Color is used for the text in Inkscape, but the package 'color.sty' is not loaded}%
    \renewcommand\color[2][]{}%
  }%
  \providecommand\transparent[1]{%
    \errmessage{(Inkscape) Transparency is used (non-zero) for the text in Inkscape, but the package 'transparent.sty' is not loaded}%
    \renewcommand\transparent[1]{}%
  }%
  \providecommand\rotatebox[2]{#2}%
  \newcommand*\fsize{\dimexpr\f@size pt\relax}%
  \newcommand*\lineheight[1]{\fontsize{\fsize}{#1\fsize}\selectfont}%
  \ifx\svgwidth\undefined%
    \setlength{\unitlength}{850.39370079bp}%
    \ifx\svgscale\undefined%
      \relax%
    \else%
      \setlength{\unitlength}{\unitlength * \real{\svgscale}}%
    \fi%
  \else%
    \setlength{\unitlength}{\svgwidth}%
  \fi%
  \global\let\svgwidth\undefined%
  \global\let\svgscale\undefined%
  \makeatother%
  \begin{picture}(1,0.83333333)%
    \lineheight{1}%
    \setlength\tabcolsep{0pt}%
    \put(0,0){\includegraphics[width=\unitlength,page=1]{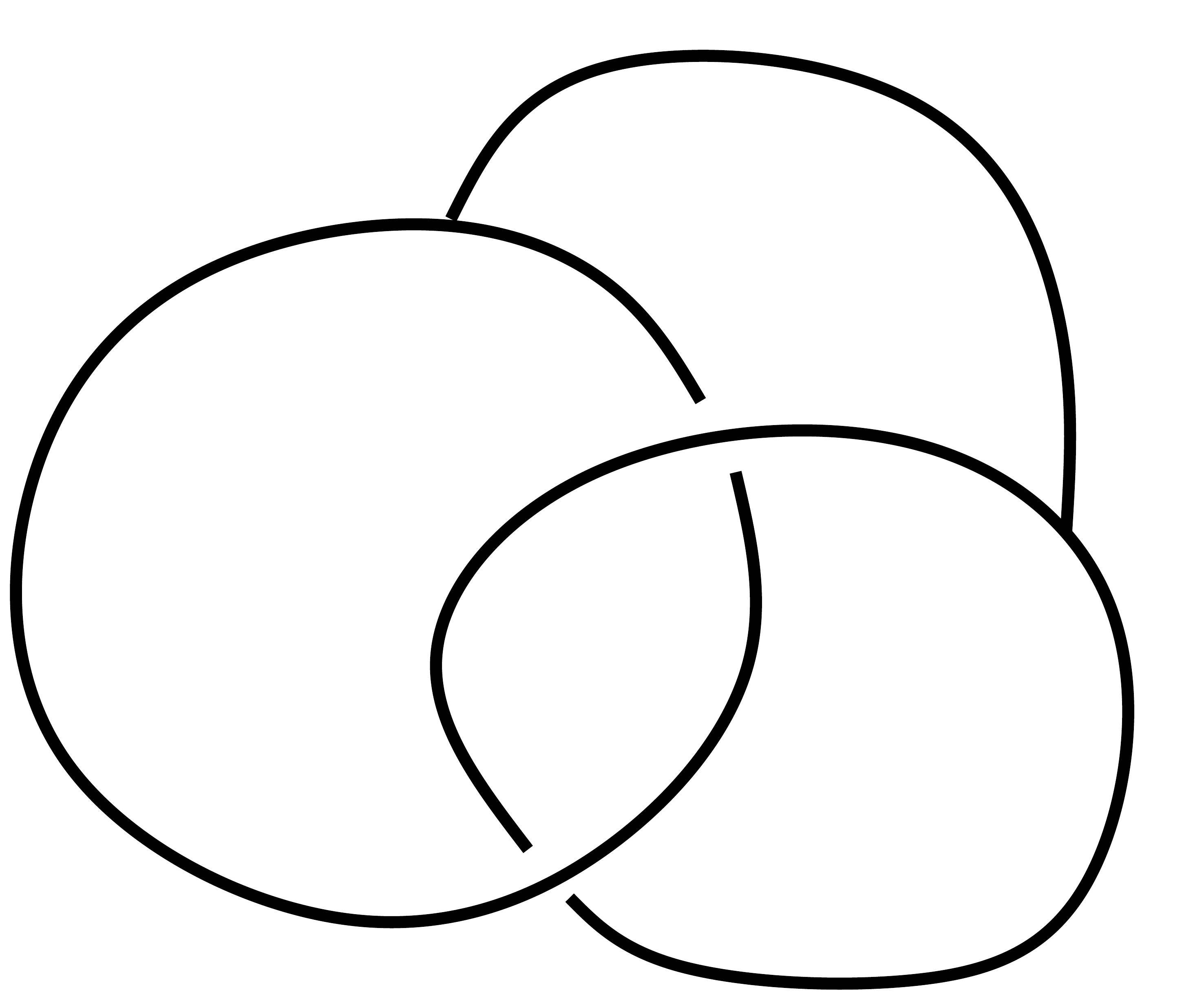}}%
    \put(0.02838688,0.04076668){\color[rgb]{0,0,0}\makebox(0,0)[lt]{\lineheight{1.25}\smash{\begin{tabular}[t]{l}{\footnotesize $L$}\end{tabular}}}}%
    \put(0.56119857,0.70318434){\color[rgb]{0,0,0}\makebox(0,0)[lt]{\lineheight{1.25}\smash{\begin{tabular}[t]{l}{\footnotesize $\gamma$}\end{tabular}}}}%
  \end{picture}%
\endgroup%

\caption{Tunnel of a Hopf link.}
\label{fig:tunnel_hopf}
\end{subfigure}
\begin{subfigure}{0.55\textwidth}
\centering
\includegraphics[scale=.12]{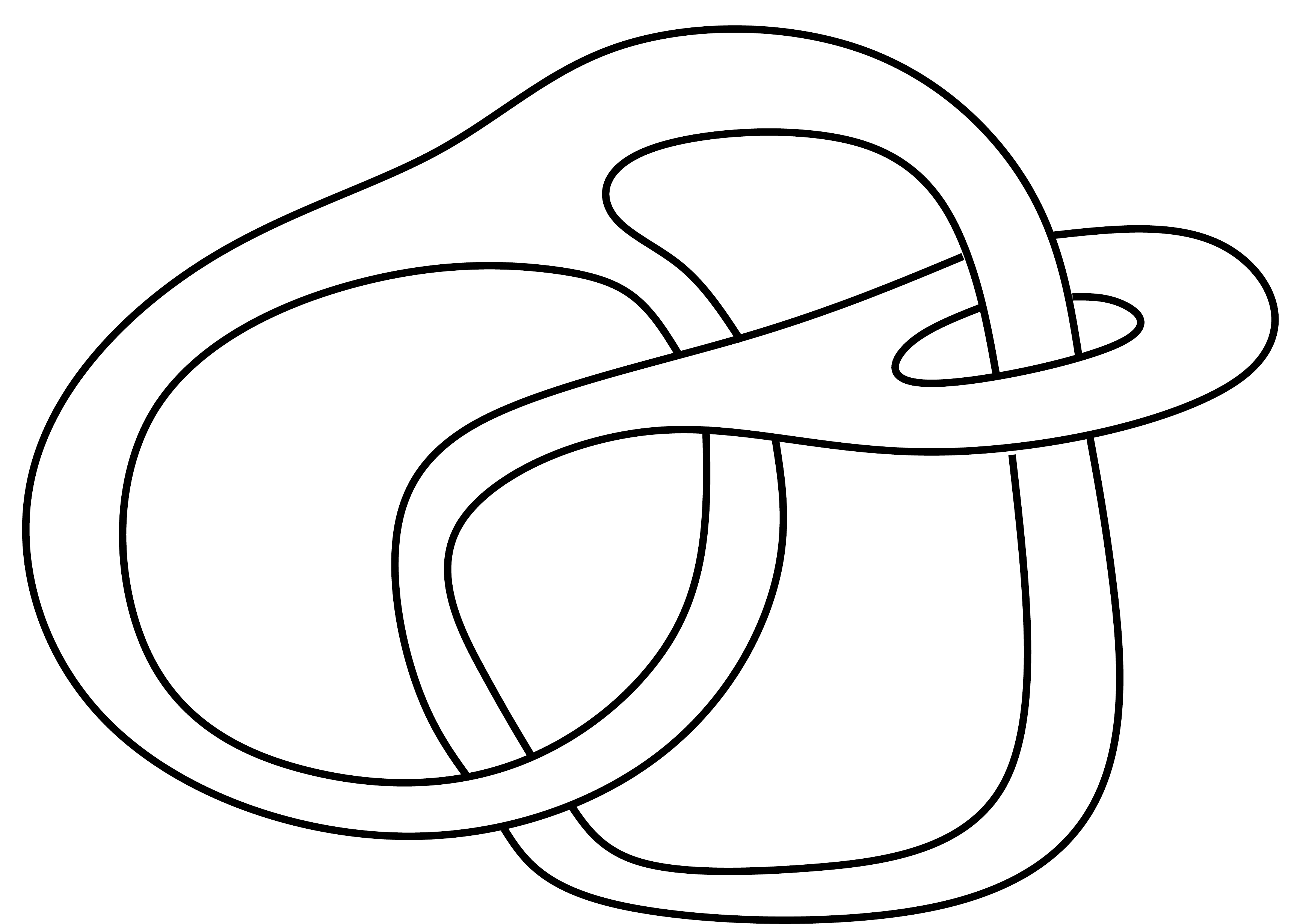} 
\caption{Ishii et al.'s $4_1$.}
\label{fig:hk4_1}
\end{subfigure}
\caption{}
\end{figure}

\smallskip

\textbf{Case $2$: $A$ is of type $2$-$2$.}
As with the previous case, we first prove 

\centerline{
{\bf Claim: $\pSym\HK=1$.}
} 

Given an element $f\in \pAut{\sphere,\HK}$,
it may be assumed $f(A)=A$. Denote by $\alpha_1, \alpha_2$ 
the components of $\partial A$ and fix an orientation of $\alpha_1,\alpha_2$.
$f$ cannot permute $\{\alpha_1,\alpha_2\}$, but it might
reverse their orientation (simultaneously) (Fig.\ \ref{fig:typetwoii}).
By Theorem \ref{teo:koda_finiteness}, 
$[f]$ is a torsion element in $\pMCG{\sphere,\HK}$, so
$[f\vert_{\partial\HK}]$ is also a torsion element 
in $\pMCG{\partial \HK}$ via the injection \eqref{eq:composition_injections}.

\textbf{Subcase $1$: $f$ preserves the orientation of $\alpha_1,\alpha_2$.}
Since the restriction of $f$   
\[f\vert_{\partial\HK}\in r^{-1}(\Emb{\alpha_1\cup\alpha_2}{\partial\HK}_0)
\subset \pAut{\partial\HK}\]
induces a torsion element in 
\[\pi_0\big(r^{-1}(\Emb{\alpha_1\cup\alpha_2}{\partial\HK}_0)\big),\] 
$[f\vert_{\partial \HK}]$ is a torsion element in 
$\pi_0(\pAut{\partial\HK,\alpha_1,\alpha_2})$ by Lemma \ref{lm:embeddings_in_surfaces}.

Now, consider the cutting homomorphism $\mathsf{cut}$
\cite[Proposition $3.20$]{FarMar:11}:
\begin{equation}\label{eq:cutting_homo}
0\rightarrow <t_1,t_2>\rightarrow \pMCG{\partial \HK, \alpha_1,\alpha_2}
\xrightarrow{\mathsf{cut}} \pMCG{\partial\HK-\alpha_1\cup\alpha_2},
\end{equation}
where the Dehn twists $t_i$ about $\alpha_i,i=1,2$ generate
its kernel.
 
Since $[f\vert_{\partial \HK-\alpha_1\cup \alpha_2}]$
is a torsion element in $\pMCG{\partial\HK-\alpha_1\cup\alpha_2}$,
and every non-trivial torsion element in $\pMCG{\partial\HK-\alpha_1\cup\alpha_2}$,
the positive mapping class group of a four-times-punctured sphere,
permutes some punctures, $[f\vert_{\partial \HK-\alpha_1\cup\alpha_2}]$ 
is the trivial element in $\pMCG{\partial\HK-\alpha_1\cup\alpha_2}$. Thus
by \eqref{eq:cutting_homo}, 
$[f\vert_{\partial\HK}]$ is trivial 
in $\pMCG{\partial\HK,\alpha_1,\alpha_2}$ as well. Thus
$f$ can be isotoped in $\pAut{\sphere,\HK,\alpha_1,\alpha_2}$
such that $f\vert_{\partial\HK}=\id$. The fact \cite{Hat:76}, \cite[Lemma $2.3$]{ChoKod:13} 
\[ \pMCG{\HK,\rel \HK}=1=\pMCG{\Compl\HK,\rel \Compl\HK}\] 
implies that $f$ can be further isotoped to $\id$ in $\pAut{\sphere,\rel \partial \HK}$, so $[f]\in \pMCG{\sphere,\HK}$ is trivial.
\cout{
\textbf{Subcase: $f$ does not reverse the orientation of $\alpha_1,\alpha_2$.}
In this case, $f\vert_{\partial_{\HK}-\alpha_1\cup\alpha_2}$ 
is necessarily trivial in $\pMCG{\partial\HK-\alpha_1\cup\alpha_2}$
since $\partial\HK-\alpha_1\cup\alpha_2$ is a four-times-punctured sphere,
and non-trivial torsion elements permute some of the punctures.
That means $f\vert_{\partial \HK}$
is trivial in $\pMCG{\partial\HK,\{\alpha_1,\alpha_2\}}$ as well,
and hence $f$ is isotopic to $\id$
in $\pSym{\HK,\alpha_1,\alpha_2}$.
}

\textbf{Subcase $2$: $f$ reverses the orientation of $\alpha_1,\alpha_2$.}
The preceding argument implies that $[f]$ is 
an element of order $2$ in 
$\pMCG{\sphere,\HK,\alpha_1,\alpha_2}$, and hence 
\[
[f\vert_{\partial\HK-\alpha_1\cup\alpha_2}]
\in\pMCG{\partial\HK-\alpha_1\cup\alpha_2}
\]
is an even element of order two.
By Lemma \ref{lm:even_elements_order_two}, one can isotope
$f\vert_{\partial\HK-\alpha_1\cup\alpha_2}$ such that
it preserves a loop $\alpha_3$ separating punctures induced
by $\alpha_1$ and punctures induced by $\alpha_2$.
Following from the isotopy extension
theorem \cite[Theorem $2$]{HudZee:64}, 
it may be assumed that $f\vert_{\partial\HK}$
preserves $\alpha_3\subset \partial\HK$, 
which separates $\alpha_1$ and $\alpha_2$ 
and therefore
bounds a separating meridian disk $D_3$
as $\alpha_2$ is the boundary of a 
meridian disk $D_2$ (Fig.\ \ref{fig:typetwoii}).

%
%
%
Let $D_1$ be a meridian disk of $\HK$
disjoint from $D_3$ with $D_1\cap \alpha_1$ a point. 
Then $f$ can be isotoped in $\pAut{\sphere,\HK,\alpha_1,\alpha_2}$ 
such that $f$ preserves $D_i$, $i=1,2,3$.
Let $D_2^\pm$ be non-parallel meridian disks of $\HK_A$ 
that are parallel to $D_2$ in $\HK$. 
Then since it may be assumed that $f(\rnbhd{A})=\rnbhd{A}$ by
the uniqueness of $A$, one can view $f$ 
as an element in $\pAut{\sphere,\HK_A,D_2^+\cup D_2^-\cup D_3}$.
Furthermore, $f(D_2^\pm)=D_2^\mp$ because 
$f$ reverses the orientation of $\alpha_2$.

Denote by $\Gamma_A$ the dual $\theta$-graph 
induced by $(\HK_A,D_2^+\cup D_2^-\cup D_3)$, and
by $K$ the constituent knot dual to $D_2^\pm$, 
and by $\gamma$ the arc dual to $D_3$.
Then $\gamma$ is a tunnel of $K$.
Via the isomorphism 
\[\pSym{\Gamma_A}\simeq \pSym{\HK_A,D_2^+\cup D_2^-\cup D_3},\]
$f$ can be isotoped such that
it preserves $K,\gamma$ and swaps
two arcs of $K$. Applying Lemma \ref{lm:asymmetry_unknotted_spatial_graph},
we see $K$ is a $2$-bridge knot and $\gamma$ is the lower or upper tunnel,
and by the classification of tunnels of a $2$-bridge knot, 
$\HK$ is a regular neighborhood of 
the union of $K$ and one of 
the four other tunnels \cite{Kob:99}, contradicting $\pair$ is non-trivial. 

\smallskip
\centerline{\bf Claim:
 $\Sym\HK=1.$
 }    
Let $f\in \Aut{\sphere,\HK}$ be orientation-reversing.   
It may be assumed that $f(\mathfrak{N}(A))=\mathfrak{N}(A)$,
and $f^2$ is isotopic to $\id$ in $\Aut{\sphere,\HK,\mathfrak{N}(A)}$
by Lemma \ref{lm:A_preserving_homo} and \eqref{eq:isomorphism_sym_hk_A}.

Denote by $A^\pm$ the two annular components of 
$\mathfrak{N}(A)\cap\partial\HK_A$.
Then $[f\vert_{\partial \HK_A}]$ and 
$[f\vert_{\partial\HK_A-A^+\cup A^-}]$
are of order two in $\MCG{\partial\HK,A^+\cup A^-}$ and
$\MCG{\partial\HK_A-A^+\cup A^-}$, respectively,
the latter being the mapping class group of a four-times-punctured
sphere. 
Since $f$ either swaps $A^+,A^-$ or preserves them,
$[f\vert_{\partial\HK_A-A^+\cup A^-}]$ is even.  
\cout{
, and is therefore
conjugate to 
\[
\big(
\begin{bmatrix}
-1&0\\
0&1
\end{bmatrix},
x\big)\in \PGL(2,\mathbb{Z})\ltimes \big(\mathbb{Z}_2\times\mathbb{Z}_2\big)\simeq \MCG{\partial\HK_A-A^+\cup A^-},
\] 
for some $x\in \mathbb{Z}_2\times \mathbb{Z}_2$.
}

By Lemma \ref{lm:even_elements_order_two}, $f\vert_{\partial\HK_A-A^+\cup A^-}$
can be isotoped such that it preserves a loop $\alpha_0$
that separates one boundary component of $A^+$ (resp.\ of $A^-$)
from the other. By the isotopy extension theorem 
\cite[Theorem $2$]{HudZee:64},  
$f\vert_{\partial\HK_A}$ can be isotoped
in $\Aut{\partial\HK_A, A^+\cup A^-}$, with $f\vert_{A^+\cup A^-}$ fixed, 
such that $f(\alpha_0)=\alpha_0$. 

Let $D_2^\pm\subset \HK_A$ 
be non-separating disks bounded by the core of $A^\pm$,
and $D_0$ a non-separating disk bounded by $\alpha_0$.
$\{D_2^+,D_2^-,D_0\}$ gives a meridian disk system of $\HK_A$,
and $f$ can be further isotoped in $\Aut{\sphere,\HK_A}$ 
such that it preserves $D_2^+\cup D_2^-$ and $D_0$.
Denote by $\Gamma_A$ the spatial graph associated to 
$(\HK_A,D_2^+,D_2^-,D_0)$, and by $K$
the constituent knot dual to $D_2^\pm$, and by $\gamma$ the arc
dual to $D_0$, which is a tunnel of $K$.

By Lemma \ref{lm:symmetry_groups_hk_gamma}, we have the isomorphism
\[\Sym{\Gamma_A}\simeq \Sym{\HK_A,D_2^+\cup D_2^-\cup D_3},\]
and hence $f$ induces an orientation-reversing homeomorphism
preserving $K,\gamma$. However, this implies $\gamma$ 
is the tunnel of a trivial knot $K$ by Lemma \ref{lm:asymmetry_unknotted_spatial_graph}.
In particular, $D_2^\pm$ are primitive disks of the Heegaard splitting
induced by $\partial \HK_A\subset\sphere$; 
therefore $A$ is $\partial$-compressible---or, in fact, $\pair$ is trivial.
Thus $\Sym\HK=1$. 
\end{proof}

\end{document}